\documentclass[a4paper,12pt]{amsart}
\usepackage{amsfonts}
\usepackage{amssymb}
\usepackage{ifthen}
\usepackage{graphicx}
\usepackage[usenames]{color}

\nonstopmode \numberwithin{equation}{section}
\newtheorem{thm}{Theorem}[section]

\newtheorem{cor}{Corollary}[section]
\newtheorem{lem}{Lemma}[section]
\newtheorem{prop}{Proposition}[section]

\newtheorem{claim}{Claim}
\newtheorem{conj}[equation]{Conjecture}

\theoremstyle{definition}
\newtheorem{defn}{Definition}[section]
\newtheorem{examp}{Example}[section]
\newtheorem{prob}[equation]{Problem}
\newtheorem{ques}{Question}[section]
\newtheorem{rem}{Remark}[section]

\newcounter {own}
\def\theown {\thesection       .\arabic{own}}

\newenvironment{pf}[1][]{%
 \vskip 3mm
 \noindent
 \ifthenelse{\equal{#1}{}}%
  {{\slshape Proof. }}%
  {{\slshape #1.} }%
 }%
{\qed\bigskip}

\newcounter{alphabet}

\newenvironment{Thm}[1][]{\refstepcounter{alphabet}%
\bigskip%
\noindent%
{\bf Theorem \Alph{alphabet}}%
\ifthenelse{\equal{#1}{}}{}{ (#1)}%
{\bf .} \itshape}{\vskip 8pt}

\newcounter{alphabet2}

\newenvironment{Lem}[1][]{\refstepcounter{alphabet}%
\bigskip%
\noindent%
{\bf Lemma \Alph{alphabet}}%
{\bf .} \itshape}{\vskip 8pt}

\newcommand{\IC}{{\mathbb C}}

\newcommand{\ID}{{\mathbb D}}

\makeatletter
\@namedef{subjclassname@2020}{%
  \textup{2020} Mathematics Subject Classification}
\makeatother

\def\be{\begin{equation}}
\def\ee{\end{equation}}

\newcommand{\bee}{\begin{enumerate}}
\newcommand{\eee}{\end{enumerate}}

\newcommand{\blem}{\begin{lem}}
\newcommand{\elem}{\end{lem}}
\newcommand{\bthm}{\begin{thm}}
\newcommand{\ethm}{\end{thm}}
\newcommand{\bcor}{\begin{cor}}
\newcommand{\ecor}{\end{cor}}
\newcommand{\beg}{\begin{examp}}
\newcommand{\eeg}{\end{examp}}
\newcommand{\begs}{\begin{examples}}
\newcommand{\eegs}{\end{examples}}
\newcommand{\bdefe}{\begin{defn}}
\newcommand{\edefe}{\end{defn}}
\newcommand{\bprob}{\begin{prob}}
\newcommand{\eprob}{\end{prob}}
\newcommand{\bques}{\begin{ques}}
\newcommand{\eques}{\end{ques}}
\newcommand{\bei}{\begin{itemize}}
\newcommand{\eei}{\end{itemize}}

\newcommand{\bca}{\begin{case}}
\newcommand{\eca}{\end{case}}
\newcommand{\bcl}{\begin{claim}}
\newcommand{\ecl}{\end{claim}}

\newcommand{\bcon}{\begin{conj}}
\newcommand{\econ}{\end{conj}}
\newcommand{\bcons}{\begin{conjs}}
\newcommand{\econs}{\end{conjs}}
\newcommand{\bprop}{\begin{prop}}
\newcommand{\eprop}{\end{prop}}
\newcommand{\br}{\begin{rem}}
\newcommand{\er}{\end{rem}}
\newcommand{\brs}{\begin{rems}}
\newcommand{\ers}{\end{rems}}
\newcommand{\bo}{\begin{obser}}
\newcommand{\eo}{\end{obser}}
\newcommand{\bos}{\begin{obsers}}
\newcommand{\eos}{\end{obsers}}
\newcommand{\bpf}{\begin{pf}}
\newcommand{\epf}{\end{pf}}
\newcommand{\ba}{\begin{array}}
\newcommand{\ea}{\end{array}}
\newcommand{\beq}{\begin{eqnarray}}
\newcommand{\beqq}{\begin{eqnarray*}}
\newcommand{\eeq}{\end{eqnarray}}
\newcommand{\eeqq}{\end{eqnarray*}}

\begin{document}

\bibliographystyle{amsplain}

\title[Sharp bounds for the growth and distortion]
{Sharp bounds for the growth and distortion of the analytic part of convex $K$-quasiconformal harmonic mappings}
\author{Peijin Li}
\address{Peijin Li,
 Department of Mathematics,
Hunan First Normal University, Changsha, Hunan 410205, People's Republic of China}
\email{wokeyi99@163.com}

\author{Saminathan Ponnusamy
}
\address{S. Ponnusamy, Department of Mathematics,
Indian Institute of Technology Madras, Chennai-600 036, India. }


\email{samy@iitm.ac.in}

\date{\today}

\subjclass[2020]{Primary 31A05; Secondary 30C55, 30C62}
\keywords{$K$-quasiconformal harmonic mappings; growth and distortion; Univalent functions; Convex functions.}


\begin{abstract}
The main aim of this paper is to obtain the sharp upper and lower bounds for the growth and distortion of the analytic part $h$ of sense-preserving convex $K$-quasiconformal harmonic mappings.
\end{abstract}


\maketitle 

\section{Introduction and Preliminaries}\label{csw-sec1}


A complex-valued function $f$ defined on  the unit disk $\mathbb{D}:=\{z\in \mathbb{C}:\, |z|< 1\}$ is called harmonic if it satisfies the
Laplace equation $\Delta f =0$, where $\Delta$ denotes the Laplacian operator
$\Delta  = 4 f_{z\overline{z}} = f_{xx} + f_{yy}.
$
If $f$ is normalized by the conditions $f(0)=0$ and $f_{z}(0)=1$, then $f$ has a canonical representation $f = h + \overline{g}$, where,
\be\label{111}
h(z)=z+\sum_{k=2}^{\infty}a_kz^k~\mbox{ and }~ g(z)=\sum_{k=1}^{\infty}b_kz^k,
\ee
and we denote the class of all such functions by $\mathcal{H}$. A locally univalent function $f= h + \overline{g}\in \mathcal{H}$ is called sense-preserving if the Jacobian $J_f(z)=|h'(z)|^2-|g'(z)|^2$ is positive in $\ID$.

\subsection{Growth and distortion of convex analytic functions}
Let ${\mathcal C}$ be the class of analytic functions $h$ in $\ID$ such that $h(0)=0$, $h'(0)=1$,   $h$ is univalent in $\ID$ and the range $h(\ID)$ is a convex domain.
It is well known (see \cite[Theorem 2.2.8]{GK}, for instance) that the growth and distortion of a given $h\in{\mathcal C}$ are controlled, respectively, by the sharp estimates
\be\label{gd1}
\frac{|z|}{1+|z|}\leq|h(z)|\leq\frac{|z|}{1-|z|},\;\;z\in\ID
\ee
and
\beqq
\frac{1}{(1+|z|)^2}\leq|h'(z)|\leq\frac{1}{(1-|z|)^2},\;\;z\in\ID.
\eeqq
If we let $|z|\rightarrow1$ in the lower bound in \eqref{gd1}, we conclude that the range $h(\ID)$ of any function $h\in{\mathcal C}$ contains the disk $\ID_{1/2}$,
and the value 1/2 is sharp as the  function $l(z)$ defined by
\be\label{lz}
l(z)=\frac{z}{1-z},
\ee
or its rotations $l_{\lambda}(z)=\overline{\lambda}l(\lambda z)$ ($|\lambda|=1$), shows. Here $\mathbb{D}_r=\{z\in \mathbb{C}:\, |z|<r\}$, and note that the function $l(z)$ maps the unit disk
conformally onto the half-plane $\mathbb{H}=\{w\in\mathbb{C}: \, {\rm Re}\,w>\frac{1}{2}\}.$




\subsection{Growth and distortion of convex harmonic functions}
Let ${\mathcal S}_H$ be the class of harmonic functions $f\in\mathcal{H}$ than are univalent and sense-preserving in $\ID$.
We denote by ${\mathcal S}_H^0$,  the subclass of ${\mathcal S}_H$ with the additional condition $f_{\overline{z}}(0)=0$. Several properties of this class together with its various geometric subclasses were investigated first by Clunie and Sheil-Small \cite{CSS} and later by many researchers (cf. \cite{DuR}).

A function $f \in \mathcal{S}_H$ is said to be convex in $\ID$ if the range $f(\ID)$ is
convex, see \cite{CSS,DuR,SS}.
The class of all convex functions $f \in \mathcal{S}_H$  is denoted by ${\mathcal C}_H$.
One can similarly define ${\mathcal C}_H^0$ in the usual manner (with the assumption that $b_1=g'(0)=f_{\overline{z}}(0)=0$).


The class ${\mathcal C}_H^0$ is much wider that its analytic counterpart ${\mathcal C}$ and there is a number of difficult problems that remain unresolved in relation with functions in this class. For instance, no sharp bounds for the growth of a given $f\in{\mathcal C}_H^0$ are known.
If $f\in{\mathcal C}_H^0$, then (see \cite[p. 100]{DuR})
\be\label{gd3}
\frac{|z|}{(1+|z|)^2}\leq|f(z)|\leq\frac{|z|}{(1-|z|)^2},\;\;z\in\ID.
\ee
However, as pointed out in \cite[p. 100]{DuR}, the lower bound in \eqref{gd3} is clearly not sharp, since the disk $|w|<\frac{1}{2}$ is contained in the range of each function
of class ${\mathcal C}_H^0$   as with the case of analytic convex mappings.
The upper bound is known to be correct with respect to the order of growth, since
\beqq
\limsup _{r \rightarrow 1}(1-r)^2 \sup _{|z|=r}|L(z)| \geq \frac{1}{2},
\eeqq
where $L=H+\overline{G}\in{\mathcal C}_H^0$ is the {\it half-plane harmonic mapping} defined by
\be\label{Lz}
L(z)=\frac{1}{2}\left[l(z)+k(z)+\overline{l(z)-k(z)}\right],\;\;z\in\ID,
\ee
$l$ is given by \eqref{lz} and $k$ is the {\it Koebe function} $k(z)=\frac{z}{(1-z)^2}$, $z\in\ID$. Note also that $L(z)$ maps the unit disk
univalently onto the half-plane  $\mathbb{H}$.

In \cite[p. 100]{DuR}, the author states that it seems likely that the sharp bounds for the growth of functions in the class ${\mathcal C}_H^0$
 are attained by the function $L$. This is a problem that is still open. Recently, in \cite{martin}, Martin obtained the sharp bounds for the growth and distortion of the analytic part $h$ of functions $f=h+\overline{g}$ in the class ${\mathcal C}_H^0$.

\begin{Thm}\label{thmA}\cite{martin}
Let $f=h+\overline{g}\in{\mathcal C}_H^0$. Then, the growth and distortion of
the analytic part $h$ of $f$ are subject, respectively, to the bounds
\beqq
\frac{2|z|+|z|^2}{2(1+|z|)^2} \leq|h(z)| \leq \frac{2|z|-|z|^2}{2(1-|z|)^2}, \;\; z \in \mathbb{D},
\eeqq
and
\beqq
\frac{1}{(1+|z|)^3} \leq\left|h^{\prime}(z)\right| \leq \frac{1}{(1-|z|)^3}, \;\; z \in \mathbb{D}.
\eeqq
The estimates are sharp: equality holds in any of the inequalities if and
only if $f$ equals the harmonic half-plane mapping \eqref{Lz} or some of its
rotations $L_{\lambda}$ given by $L_{\lambda}(z)=\overline{\lambda}L(\lambda z)$, where $|\lambda|=1$. That is, if and
only if $h$ equals the analytic part $H$ of $L$ or some of its rotations.
\end{Thm}

A natural question is the following:

\bques\label{ques1}
What is counterpart of Theorem~A 
for the class ${\mathcal C}_H^0(K)$ of convex $K$-quasiconformal harmonic mappings?
\eques

The definition of ${\mathcal C}_H^0(K)$ is given below. 
The aim of this article is to present an affirmative answer to this question.

\section{Main result and discussion}

If a sense-preserving and locally univalent harmonic mapping $f=h+\overline{g}$ on $\ID$ satisfies the condition
$ |\omega_f(z)| \leq k<1$ for $\ID$,  then $f$ is called $K$-quasiregular harmonic mapping in $\ID$, where $\omega_f =g'/h'$ and
$K=\frac{1+k}{1-k}\geq 1$. We say that $f$ belongs to the class ${\mathcal S}_H(K)$ of $K$-quasiconformal harmonic mappings, if $f\in {\mathcal S}_H$ and $K$-quasiregular harmonic mapping  in $\ID$ (cf. \cite{CP-22}).
 As with the standard practice, a function $f$ is called a $K$-quasiconformal harmonic mapping if it belongs to
${\mathcal S}_H(K)$ for some $K\geq 1$. Also, we define
$${\mathcal S}_H^0(K):={\mathcal S}_H(K)\cap {\mathcal S}_H^0.
$$
It is useful to introduce the following analogous notations:
$$
{\mathcal C}_H(K):={\mathcal S}_H(K)\cap {\mathcal C}_H ~\mbox{ and }~{\mathcal C}_H^0(K):={\mathcal S}_H(K)\cap {\mathcal C}_H^0.
$$

In \cite{LP}, the authors investigated the classes ${\mathcal C}_H(K)$ and ${\mathcal C}_H^0(K)$,  and also introduced {\it half-plane $K$-quasiconformal harmonic mapping} $L_k(z)$ defined by
\be\label{Lkz}
L_k(z)=H_k(z)+\overline{G_k(z)},
\ee
where
\be\label{hkz}
H_k(z)=\frac{z}{(1-k)(1-z)}+\frac{k}{(1-k)^2} \log \left(\frac{1-z}{1-k z}\right),
\ee
and
\beqq
G_k(z)=\frac{-k z}{(1-k)(1-z)}-\frac{k}{(1-k)^2} \log \left(\frac{1-z}{1-k z}\right),
\eeqq
and showed that $L_k\in{\mathcal C}_H^0(K)$. Given $\lambda\in \IC$ with $|\lambda|=1$, the {\it rotation} $L_k^{\lambda}$ of $L_k$ is $L_k^{\lambda}(z)=\overline{\lambda}L_k(\lambda z)$. It is easy to check that $L_k^{\lambda}\in{\mathcal C}_H^0(K)$ for all $|\lambda|=1$. 
The {\it rotation} $H_k^{\lambda}$ of $H_k$ is defined by
\be\label{hkz0}
H_k^{\lambda}(z)=\overline{\lambda}H_k(\lambda z)
=\frac{\overline{\lambda}\lambda z}{(1-k)(1-\lambda z)}+\frac{\overline{\lambda}k}{(1-k)^2} \log \left(\frac{1-\lambda z}{1-k\lambda z}\right),
\ee
for $\in \mathbb{D}$, and
\be\label{hkz00}
\left(H_k^{\lambda}\right)'(z)=\frac{1}{(1-\lambda z)^2(1-k\lambda z)}.
\ee

If $\lambda=\lambda_1=\frac{\overline{z}}{|z|}$  $(z\neq 0)$ in \eqref{hkz0} and \eqref{hkz00} respectively, we have
\be\label{hkz1}
H_k^{\lambda_1}(z)
=\overline{\lambda_1}\left[\frac{z}{(1-k)(1-|z|)}+\frac{k}{(1-k)^2} \log \left(\frac{1-|z|}{1-k |z|}\right)\right]
\ee
and
$$ \left(H_k^{\lambda_1}\right)'(z)=\frac{1}{(1-|z|)^2(1-k|z|)}.$$

If $\lambda=\lambda_2=-\frac{\overline{z}}{|z|}$ $(z\neq 0)$ in \eqref{hkz0} and \eqref{hkz00} respectively, we have
\be\label{hkz2}
H_k^{\lambda_2}(z)
=-\overline{\lambda_2}\left[\frac{z}{(1-k)(1+|z|)}+\frac{k}{(1-k)^2} \log \left(\frac{1+k |z|}{1+|z|}\right)\right]
\ee
and
$$ \left(H_k^{\lambda_2}\right)'(z)=\frac{1}{(1+|z|)^2(1+k|z|)}.$$

In this paper, we consider the sharp upper and lower bounds for the growth and distortion of the analytic part $h$ of
$f=h+\overline{g}\in{\mathcal C}_H^0(K)$. The following theorem answers Question \ref{ques1} affirmatively.

\begin{thm}\label{thm1}
Let $f=h+\overline{g}\in{\mathcal C}_H^0(K)$. Then, for $z \in \mathbb{D}$, the growth and distortion of
the analytic part $h$ of $f$ are subject, respectively, to the bounds
\be\label{gd4}
B(k, |z|)\leq|h(z)| \leq A(k, |z|)
\ee
and
\be\label{gd5}
\frac{1}{(1+|z|)^2(1+k|z|)} \leq\left|h^{\prime}(z)\right| \leq \frac{1}{(1-|z|)^2(1-k|z|)}.
\ee
Here
$$A(k, |z|)=\frac{|z|}{(1-k)(1-|z|)}+\frac{k}{(1-k)^2}\log\left(\frac{1-|z|}{1-k|z|}\right)$$
and
$$B(k, |z|)=\frac{|z|}{(1-k)(1+|z|)}+\frac{k}{(1-k)^2}\log\left(\frac{1+k|z|}{1+|z|}\right).$$
The estimates are sharp: the upper bound equality holds in any of the inequalities if and
only if $h$ equals $H_k^{\lambda_1}(z)$ defined by \eqref{hkz1} and the lower bound equality holds in any of the inequalities if and
only if $h$ equals $H_k^{\lambda_2}(z)$ defined by \eqref{hkz2}.
\end{thm}

\begin{rem}\label{rem-3}
If $k\rightarrow1^-$, then Theorem \ref{thm1} shows that
$$\lim_{k\rightarrow1^-}A(k, |z|)=\frac{2|z|-|z|^2}{2(1-|z|)^2}
~\mbox{ and }~
\lim_{k\rightarrow1^-}B(k, |z|)=\frac{2|z|+|z|^2}{2(1+|z|)^2}.$$
In this limiting case, Theorem \ref{thm1} clearly coincides with Theorem~A. 
\end{rem}

We  present the proof of  Theorem \ref{thm1} in Section \ref{sec-2} and, in order to do this, we need the following result.

\begin{Lem}$($\cite[p.~51]{DuR}$)$\label{lemA}
If $f=h+\overline{g}\in {\mathcal C}_H$, then there exist angles $\alpha$ and $\beta$ such that
$${\rm Re}\, \{(e^{i\alpha}h'(z)+e^{-i\alpha}g'(z))(e^{i\beta}-e^{-i\beta}z^2)\}>0 ~\mbox{ for all $z\in\ID$.}
$$
\end{Lem}

\section{Proof of Theorem \ref{thm1}}\label{sec-2}


Let $f=h+\overline{g}\in{\mathcal C}_H^0(K)$ have dilatation $\omega(z)=g'(z)/h'(z)$, and let $a\in\ID$.
We know that (see \cite{martin} or \cite[P. 79]{DuR})
$$F=H+\overline{G}=\frac{f(\varphi(z))-f(a)}{h'(a)(|a|^2-1)}\in{\mathcal C}_H,$$
where
$$\varphi(z)=\frac{a-z}{1-\overline{a}z},\;\;z\in\ID.$$
By calculation, we get
$$\omega_F(z)=\frac{G'(z)}{H'(z)}=\frac{h'(a)\omega(\varphi(z))}{\overline{h'(a)}},$$
and thus, $|\omega_F(z)|=|\omega(\varphi(z))|\leq k$ for $z\in \ID$. Therefore, $F\in{\mathcal C}_H(K)$.

It is known from the proof of Corollary 1.3 in \cite{LP} that each $F=H+\overline{G} \in {\mathcal C}_H(K)$ has the form $F=F_0+\overline{G'(0)F_0}$ for some $F_0=H_0+\overline{G_0}\in {\mathcal C}_H^0(K_0)$, where $K\geq K_0$. Hence,
\be\label{F0}
F_0(z)=\frac{F(z)-\overline{\omega_F(0)F(z)}}{1-|\omega_F(0)|^2}\in{\mathcal C}_H^0(K_0),
\ee
since $G'(0)=\omega_F(0)$.
A straightforward calculation shows that
\be\label{H0}
H_0(z)=\frac{h(\varphi(z))-h(a)-\overline{\omega(a)}\big(g(\varphi(z))-g(a)\big)}{h'(a)(|a|^2-1)(1-|\omega(a)|^2)}
\ee
and
\beqq
G_0(z)=\frac{g(\varphi(z))-g(a)-\omega(a)\big(h(\varphi(z))-h(a)\big)}{\overline{h'(a)}(|a|^2-1)(1-|\omega(a)|^2)}.
\eeqq
Then the dilatation $\omega_{F_0}$ of $F_0$ is
\be\label{WF0}
\omega_{F_0}(z)=\frac{h'(a)}{\overline{h'(a)}}\frac{\omega(\varphi(z))-\omega(a)}{1-\overline{\omega(a)}\omega(\varphi(z))}.
\ee

\subsection{Proof of  Theorem \ref{thm1}}
Let $f=h+\overline{g}\in{\mathcal C}_H^0(K)$ have dilatation $\omega(z)=g'(z)/h'(z)$.
We recall  that the dilatation $\omega(z)=g'(z)/h'(z)$ satisfies $\omega(0)=0$ and $|\omega(z)|\leq k$ and
thus, by the classical Schwarz Lemma,  we have
$|\omega(z)|\leq k|z|$ for $z\in \ID$ and $|\omega'(0)|\leq k$.

We first prove the upper bound in \eqref{gd5}. Appealing now to Lemma~B,
we conclude that there exist angles $\alpha$ and $\beta$ such that
$${\rm Re}\,\{(e^{i\alpha}h'(z)+e^{-i\alpha}g'(z))(e^{i\beta}-e^{-i\beta}z^2)\}>0, ~\mbox{ for $z\in\ID$.}~
$$
The last relations is equivalent to
$${\rm Re}\,\{e^{i\gamma} p(z) \}>0, ~\mbox{ i.e. }~ e^{i\gamma} p(z)\prec \frac{e^{i\gamma} +e^{-i\gamma} z}{1-z}, \mbox{ for $z\in\ID$.}~
$$
where $p(z)=(h'(z)+e^{-2i\alpha}g'(z))(1-e^{-2i\beta}z^2)$, $\gamma =\alpha +\beta$ and $\prec$ denotes the usual subordination (cf. \cite{DuR}).
By the definition of subordination, the above relation may be written equivalent as follows: there exists an analytic
function $\delta$ in the unit disk with $\delta(0)=0$ and $\delta(\ID)\subset\ID$ such that for $z\in\ID$,
$$h^{\prime}(z)+e^{-2 i \alpha} g^{\prime}(z)=\frac{1+e^{-2 i(\alpha+\beta)} \delta(z)}{1-\delta(z)} \frac{1}{1-e^{-2 i \beta} z^2}.$$
As a consequence of the triangle inequality, we have
\be\label{h}
\left|h^{\prime}(z)+e^{-2 i \alpha} g^{\prime}(z)\right|=|h^{\prime}(z)|\cdot|1+e^{-2 i \alpha}\omega(z)|
\leq \frac{1+|z|}{1-|z|} \frac{1}{1-|z|^2}=\frac{1}{(1-|z|)^2},
\ee
where $\omega(z)=g'(z)/h'(z)$. Therefore, we have
$$|h^{\prime}(z)|\leq\frac{1}{(1-|z|)^2}\frac{1}{|1+e^{-2 i \alpha}\omega(z)|}\leq\frac{1}{(1-|z|)^2(1-k|z|)}.$$
This proves the upper bound in \eqref{gd5}.

To prove that the upper bound in \eqref{gd4} is satisfied, let us fix $z=re^{ i \theta}\in\ID\setminus\{0\}$, so
that $|z|=r\in(0, 1)$. and observe that, since $h(0)=0$,
$$h(z)=\int_0^r h^{\prime}\left(\rho e^{i \theta}\right) e^{i \theta} d \rho.$$
Thus,
\beqq
|h(z)| &\leq& \int_0^r\left|h^{\prime}\left(\rho e^{i \theta}\right)\right||d \rho| \leq \int_0^r \frac{d \rho}{(1-\rho)^2(1-k\rho)}\\
&=&\frac{r}{(1-k)(1-r)}+\frac{k}{(1-k)^2}\log\left(\frac{1-r}{1-kr}\right)
\eeqq
which  proves the upper bound in \eqref{gd4}.

Let us now prove the lower bound in \eqref{gd5}. To do so, let $a\in\ID$.
Since $F_0(z)\in{\mathcal C}_H^0(K_0)$ defined by \eqref{F0},
we can argue as in \eqref{h} that there exist angles $\alpha$ such that for all $z\in\ID$,
$$\left|H_0^{\prime}(z)+e^{-2 i \alpha} G_0^{\prime}(z)\right|\leq \frac{1}{(1-|z|)^2}.$$
In particular, if we set $z = a$, we get
\be\label{H1}
\left|H_0^{\prime}(a)+e^{-2 i \alpha} G_0^{\prime}(a)\right|=|H_0^{\prime}(a)|\cdot|1+e^{-2 i \alpha}\omega_{F_0}(a)|
\leq \frac{1}{(1-|a|)^2},
\ee
where $\omega_{F_0}$ is the dilatation of $F_0$ given by \eqref{WF0}. Since $\omega(0)=0$, we have $\omega_{F_0}(a)=|\omega(a)|$.
Therefore, it follows from \eqref{H1} that
\be\label{H2}
|H_0^{\prime}(a)|(1-|\omega(a)|) \leq|H_0^{\prime}(a)|\cdot|1+e^{-2 i \alpha}\omega_{F_0}(a)|
\leq \frac{1}{(1-|a|)^2}.
\ee
Using the formula for $H_0$ in \eqref{H0} and \eqref{H2}, a straightforward calculation shows that for any $a\in\ID$,
$$|H_0^{\prime}(a)|=\frac{1}{|h'(a)|(1-|a|^2)^2(1-|\omega(a)|^2)}\leq \frac{1}{(1-|a|)^2(1-|\omega(a)|)}.$$
This is equivalent to
$$|h'(a)|\geq\frac{1}{(1+|a|)^2(1+|\omega(a)|)}\geq\frac{1}{(1+|a|)^2(1+k|a|)},$$
since $|\omega(a)|\leq k|a|$. This proves the lower bound in \eqref{gd5}.

For $r\in(0, 1)$, let $m_r=\min_{|z|=r}|h(z)|$. By the continuity of $h$ in $\ID$,
there must exists a point $z_0$ of modulus $r$ such that $|h(z_0)| = m_r$.
Let the segment which connects $0$ and $h(z_0)$ is $s$.
Since $h$ is univalent in the unit disk, we have that $C$, which is the pre-image under $h$ of the segment $s$, is an simple analytic curve which joins $0$ and $h(z_0)$. Thus,
\beqq
|h(z_0)| &=& \left|\int_C h'(\zeta)d\zeta\right|=\int_C |h'(\zeta)||d\zeta|\geq\int_0^r \frac{d \rho}{(1+\rho)^2(1+k\rho)}\\
&=&\frac{r}{(1-k)(1+r)}+\frac{k}{(1-k)^2}\log\left(\frac{1+kr}{1+r}\right).
\eeqq
This shows that if $|z| = r$, then
$$|h(z)|\geq|h(z_0)|\geq\frac{r}{(1-k)(1+r)}+\frac{k}{(1-k)^2}\log\left(\frac{1+kr}{1+r}\right).$$
This proves the lower bound in \eqref{gd4}.

The sharpness are obvious follows from \eqref{hkz1} and \eqref{hkz2}.
This completes the proof of Theorem \ref{thm1}.
\qed

\subsection*{Acknowledgments}
The research was partly supported by the Natural Science Foundation of  China (No. 12371071).

\subsection*{Conflict of Interest Statement}
The authors declare that they have no conflict of interest, regarding the publication of this paper.

\subsection*{Data Availability Statement}
The authors declare that this research is purely theoretical and does not associate with any datas.


\end{document}